\documentclass[12pt,reqno,oneside]{amsart}
\usepackage{amsmath,amsthm,amsfonts,amssymb}
\usepackage[mathscr]{eucal}
\usepackage{indentfirst}
\usepackage{url}

\theoremstyle{plain}
\newtheorem{teo}{Theorem}[section]
\newtheorem{cor}[teo]{Corollary}

\newtheorem{prop}[teo]{Proposition}

\theoremstyle{definition}
\newtheorem{defn}[teo]{Definition}
\newtheorem{exa}[teo]{Example}
\newtheorem{rem}[teo]{Remark}

\newtheorem{pro}[teo]{Proposition}

\numberwithin{equation}{section}

\def\bbR{{\mathbb R}}
\def\bbP{{\mathbb P}}
\def\bbZ{{\mathbb Z}}
\def\bbN{{\mathbb N}}
\def\bbE{{\mathbf E}}

\def\bbT{{\mathbb T}}

\def\cpp{\textit{Cone Percolation Model}}
\def\rpm{\textit{Rumor Percolation Model}}
\def\dpm{\textit{Disk Percolation Model}}
\def\hcpp{\textit{Heterogeneous Cone Percolation Process}}
\def\hofp{\textit{Homogeneous Fireworks Process}}
\def\horfp{\textit{Homogeneous Reverse Fireworks Process}}
\def\hefp{\textit{Heterogeneous Fireworks Process}}
\def\herfp{\textit{Heterogeneous Reverse Fireworks Process}}

\begin{document}
\baselineskip=22pt
\title[The rumor percolation model and its variations]{The rumor percolation model and its variations}
\author{Valdivino~V.~Junior}
\author{F\'abio~P.~Machado}
\author{Krishnamurthi~Ravishankar}
\address[F\'abio~P.~Machado]
{Institute of Mathematics and Statistics
\\ University of S\~ao Paulo \\ Rua do Mat\~ao 1010, CEP
05508-090, S\~ao Paulo, SP, Brazil - fmachado@ime.usp.br }

\noindent
\address[Valdivino~V.~Junior]
{Federal University of Goias
\\ Campus Samambaia, Goi\^ania, GO, Brazil - vvjunior@ufg.br}

\address[Krishnamurthi Ravishankar]
{NYU-ECNU Institute of Mathematical Sciences at NYU Shanghai, 3663
Zhongshan Road North, Shanghai 200062 and 1555 Century Ave, Pudong, Shanghai, China - kr26@nyu.edu}

\noindent

\thanks{Research supported by CNPq (310829/2014-3), FAPESP (09/52379-8), PNPD-Capes 536114 and Simons Foundation Collaboration grant 281207}



\date{\today}

\begin{abstract}
The study of rumor models from a percolation theory point of view has gained a few adepts in the last few years. The persistence of a rumor, which may consistently spread out throughout a population can be associated to the existence of a giant component containing the origin of a graph. That is one of the main interest in percolation theory. In this paper we present a quick review of recent results on rumor models of this type.
\end{abstract}

\maketitle

\tableofcontents
%\cleardoublepage

\section{Introduction and basic definitions}
\label{S: Introduction}

We are interested in a long-range percolation model on infinite graphs which we call the \rpm. 
Such models have recently been studied by a few authors in a series of papers.
The dynamics of the model describes the spreading of a rumor on a graph in the following way. We assign independent random \textit{radius of influence} $R_v$ to each vertex $v$ of an infinite, locally finite, connected graph $\mathcal G$. Then
we define a chain reaction on $\mathcal G$ according to the following simple rules: (1) at time zero, only the root (a fixed vertex of $\mathcal G$)
hears the rumor, (2) at time $n \ge 1$, a new vertex hears the rumor if it is a distance at most $R_v$ of some vertex $v$ that
previously heard the rumor.
We point out that similar models, are of interest in Computer Science, in particular in the area of distributed networks. One of the problems of interest is the broadcasting problem where one node has some information which it wants to pass on to other nodes. Questions of optimal algorithm for achieving this goal are of interest. This question was considered for the case where the nodes are uniformly randomly distributed on an interval $[0,L]$ and the nodes had a transmission radius of one. In \cite{RS94} asymptotically (in $L$) optimal algorithm was obtained.

\begin{defn}
\label{D: DefinicaoCPP}
The \rpm\ on $\mathcal G$.\\
Let $\mathcal G = ({\mathcal V},{\mathcal E})$ be an infinite,
locally finite, connected graph and
let $\{R_v\}_{\{ v \in {\mathcal V} \}}$ be a set of independent and identically distributed random variables.
Furthermore, for each $u \in {\mathcal V}$, we define the random sets
\begin{equation}
\label{E: defBu}
B_u = \{v \in {\mathcal V}: d(u,v) \leq R_u\}.
\end{equation}
or
\begin{equation}
\label{E: defBuC}
B_u = \{v \in {\mathcal V}, u \leq v : d(u,v) \leq R_u\}.
\end{equation} 
\noindent
With these sets we define the \rpm\ on ${\mathcal G}$, the non-decreasing sequence of random
sets $I_0 \subset I_1 \subset \cdots$
defined as $ I_0 = \{{\mathcal O}\} $ and inductively $I_{n+1} = \bigcup_{u \in I_n} B_u$
for all $ n \geq 0.$
\end{defn}

\begin{defn}
\label{D: Survival}
The \rpm\ $survival$.\\
Consider $ I = \bigcup_{n \geq 0} I_n$ be the connected component of the
origin of $\mathcal G$. Under the rumor process interpretation, $I$ is the set of vertices
which heard the rumor. We say that the process \emph{survives} (\emph{dies out}) if $|I|=\infty$ ($|I|<\infty$),
referring to the surviving event as $V.$
\end{defn}

In section 2 we review the paper of Athreya~\textit{et al}~\cite{ARS}. Instead of considering a graph
structure they consider a homogeneous Poisson point process on ${{\bbR}^d}$ and ${{\bbR}^d}^+$ with $\{R_v\}$, the \textit{box of influence}, starting
from every point $v$ of the point process in the sense of~(\ref{E: defBuC}).
They work with the concept of the coverage of
a set $(t, \infty)^d$ for some $t > 0$, the \textit{eventual coverage}.
In section 3 we review the paper of Lebenstayn and Rodriguez~\cite{LR} where authors consider the $\dpm$. While the set of \textit{radius of influence}, 
$\{R_v\}_{\{ v \in {\mathcal V} \}}$, has a geometric distribution, the graph ${\mathcal G}$ is quite general.
In their version the \textit{radius of influence} of a vertex $v \in {\mathcal G}$ goes in every possible direction
as in~(\ref{E: defBu}). In section 4 we review the papers of Junior \textit{et al}~\cite{Junior} and 
Gallo \textit{et al}~\cite{GGJR}. They work with a processes that they made known as Fireworks on $\bbN$
(direct and reverse). They studied an homogeneous version, where there is one informant per vertex
and the \textit{radius of influence} are independent and have the same distribution, and a heterogeneous version, where one of these conditions fail. In their models the \textit{radius of influence} goes like in~(\ref{E: defBuC}). In section 5 the papers of Junior {\it et al}~\cite{Junior},~\cite{Junior2} and~\cite{Junior3} are briefly reviewed. They work with the Cone Percolation model, a Fireworks model in a tree (homogeneous, spherically symmetric, periodic or Galton Watson). In all these models the the \textit{radius of influence} goes like in~(\ref{E: defBuC}).  In section 6 we review the paper of Bertacchi and Zucca~\cite{BZ}. They consider 
a type of random environment in the sense that the number of informants in each vertex of are
random.

\section{Random sets on ${{\bbR}^d}$ and ${{\bbR}^d}^+$}
\label{S: RS}

The theory of coverage processes was introduced by P. Hall~\cite{Hall} in 1988. He developed a class of stochastic processes intended to
be used as a model for \textit{binary images}, that is, images which partition ${\mathbb R}^d$ into two regions, ${\mathcal C}$ and its complement, representing the ``black'' and ``white'' parts of an image. In its basic version the process consists of a point process $P = \{\xi_1, \xi_2, \dots \}$ and a collection of random sets $\{S_1, S_2, \dots \}$. The ``black'' region ${\mathcal C}$ is then defined to be ${\mathcal C} = \cup_{i=1}^{\infty} (\xi_i + S_i).$ P. Hall~\cite{Hall} developed probabilistic results on geometrical
properties of ${\mathcal C}$, such as the size-distribution of its connected subsets. In that work the main assumptions needed to obtain explicit results is
that $P$ is an homogeneous Poisson process and the $S_i$ are independent copies of a random closed set. This version is known as the Poisson Boolean model.

%Athreya~\textit{et al}~\cite{ARS} considered a Poisson Bolean model with 
%\begin{displaymath}
%{\mathcal C} = \cup_{i=1}^{\infty} (\xi_i + {[0, \rho_i]}^d)
%\end{displaymath}
%where $\{\rho_1, \rho_2, \dots \}$
%are independent identically distributed positive random variables independent of the point process. 
%%They prove that the existence of $ t > 0$ such that
%%$ {(t, \infty)}^d \subset {\mathcal C}$ (\textit{eventual coverage})  depends on the behaviour of $xP(\rho > x)$ as $x \to \infty$ but also
%%depends on whether $d = 1$ or $d \ge 2$.

Athreya~\textit{et al}~\cite{ARS} considered two different models, both related to rumor percolation. For the first model, arising for genome analysis, they consider 
$\{X_i\}_{i \in \bbN}$ be a
$\{0,1\}$-valued time-ho\-mo\-ge\-neous Markov chain and  $\{\rho_i\}_{i \in \bbN}$
an independent and identically distributed sequence of random variables assuming values on $\bbN$, independent of the Markov chain. Let $S_i = [i, i +\rho_i]$
whenever $X_i=1$ ($\emptyset$ otherwise) and $C = \cup_{i=1}^{\infty} S_i$.

\begin{defn}
\label{D: EventuallyCovered}
We say that $\bbN$ is eventually covered by $C$ (or $C$ eventually covers $\bbN$)  if there exists a $ t \ge 1$ such
that $[t, \infty) \subseteq C$.
\end{defn}

\begin{teo}[Athreya \textit{et al}~\cite{ARS}]
\label{T: ARS1.1} Let $p_{ij} = \bbP(X_{n+1}=j | X_n=i)$. Assume that
$0< p_{00}, p_{10}<1,$
\begin{itemize}
\item[\textit{(i)}] If
\[ l = \liminf_{j \to \infty} j \bbP(\rho_1 > j) > 1, \]
then
\[ \bbP(C \hbox{eventually covers } \bbN) = 1\]
whenever
\[ \frac{p_{01}}{p_{10}+p_{01}} > \frac{1}{l}.\]
\item[\textit{(ii)}] If
\[ L = \limsup_{j \to \infty} j \bbP(\rho_1 > j) < \infty, \]
then
\[ \bbP(C \hbox{eventually covers } \bbN) = 0\]
whenever
\[ \frac{p_{01}}{p_{10}+p_{01}} < \frac{1}{L}.\]
\end{itemize}
\end{teo} 

Their second model aims to complement known results on complete coverage in stochastic geometry. For $B(0, \rho)$ the closed $d$-dimensional ball of radius $\rho$ centered at the origin,
some important previous results for the random covered region 
$\displaystyle{\cup_{i=1}^{\infty} (\xi_i + B(0,\rho_i))}$ are presented in the next two theorems.

\begin{teo}[Hall~\cite{Hall}]
\label{T: Hall}
For the Poisson Bolean model on $\bbR^d$ the space is fully covered by
$\cup_{i=1}^{\infty} (\xi_i + B(0,\rho_i))$ almost surely if and only if 
$\bbE(\rho^d) = \infty$.
\end{teo}

If instead of a Poisson point process one considers an arbitrary ergodic process,
there is the following result

\begin{teo}[Meester and Roy~\cite{MR}]
\label{T: MR}
For the Bolean model on $\bbR^d$ the space is fully covered by
$\cup_{i=1}^{\infty} (\xi_i + B(0,\rho_i))$ almost surely if 
$\bbE(\rho^d) = \infty$.
\end{teo}

Athreya \textit{et al}~\cite{ARS} take $\bbR^d_+$ and the random covered region
\begin{displaymath}
C = \cup_{\{i: \xi \in \bbR^d_+\}} (\xi_i +[0,\rho_i]^d).
\end{displaymath}

Guided by the fact that $C$ will never completely cover $\bbR^d_+$ because, for
any $ \epsilon >0,$ ${[0, \epsilon]}^d$ will not be covered by $C$ with positive probability,
they work with the notion of \textit{eventual coverage} for the orthant $\bbR^d_+$.

\begin{defn}
\label{D: REventuallyCovered}
We say that $\bbR^d_+$ is eventually covered by the Poisson Boolean
model if there exists a $ t \in (0, \infty) $ such
that ${[t, \infty)}^d \subseteq C$.
\end{defn}

With this notion Athreya \textit{et al}~\cite{ARS} are able to present the following result, considering a Poisson Bolean model
on $\bbR^d_+$. They show that eventual coverage depends on the growth rate of the distribution function of $\rho$ (even when $\bbE(\rho) = \infty$) as well as on whether $d=1$ or $d \ge 2$.

\begin{teo}[Athreya \textit{et al}~\cite{ARS}]
\label{T: Royd=1}
Assume $d=1$.
\begin{itemize}
\item[\textit{(i)}] If
\[ 0 < l :=\liminf_{x \to \infty} x \bbP(\rho > x) < \infty, \]
then there exists a $\lambda_0$ such that $0< \lambda_0 \le {1}/{l} < \infty$ and
\[ \bbP_{\lambda}(\bbR_+ \hbox{ is eventually covered by } C) = 
\left\{ \begin{array}{ll}
0 & \mbox{if $\lambda < \lambda_0$,} \\
1 & \mbox{if $\lambda > \lambda_0$;}
\end{array}
\right. \]
\item[\textit{(ii)}] If
\[ 0 < L :=\limsup_{x \to \infty} x \bbP(\rho > x) < \infty, \]
then there exists a $\lambda_1$ such that $0 < {1}/{L} \le  \lambda_1 < \infty$ and
\[ \bbP_{\lambda}(\bbR_+ \hbox{ is eventually covered by } C) = 
\left\{ \begin{array}{ll}
0 & \mbox{if $\lambda < \lambda_1$,} \\
1 & \mbox{if $\lambda > \lambda_1$;}
\end{array}
\right. \]

\item[\textit{(iii)}] If
\[ \lim_{x \to \infty} x \bbP( \rho > x) = \infty, \]
then for all $\lambda > 0, \bbR_+$ is eventually covered by $C$ 
($\bbP_{\lambda}-$a.s.);

\item[\textit{(iv)}] If
\[ \lim_{x \to \infty} x \bbP( \rho > x) = 0, \]
then for any $\lambda > 0, \bbR_+$ is eventually covered by $C$ 
($\bbP_{\lambda}-$a.s.).

\end{itemize}
\end{teo}

\begin{teo}[Athreya \textit{et al}~\cite{ARS}]
\label{T: Royd>1}
Let $d \ge 2$. For all $\lambda > 0$,
\begin{itemize}

\item[\textit{(i)}] If
\[ \liminf_{x \to \infty} x \bbP( \rho > x) >0, \]
then
\[ \bbP_{\lambda}(\bbR^d_+ \hbox{ is eventually covered by } C) = 1; \]

\item[\textit{(ii)}] If
\[ \lim_{x \to \infty} x \bbP( \rho > x) = 0, \]
then
\[ \bbP_{\lambda}(\bbR^d_+ \hbox{ is eventually covered by } C) = 0. \]
\end{itemize}
\end{teo}

It is interesting to observe that while $\bbE(\rho^d) = \infty$ guarantees complete coverage
of $\bbR^d$ by $C$, it is not sufficient to guarantee \textit{eventual coverage} for $\bbR^d_+$.
This is due to the fact that a boundary effect is present in the orthant $\bbR^d_+$
but absent in the whole space $\bbR^d$.

\section{Disk percolation}
\label{S: DPM}

Lebensztayn and Rodriguez studied a long-range percolation model on infinite graphs, the $\dpm$. They assign a random \textit{radius of influence} $R_v$ to each vertex $v$ of an infinite, locally finite, connected graph $G$, so that all the assigned radii are independent and identically distributed random variables with geometric distribution with parameter$(1-p)$, which means, satisfying
\[
\bbP(R = k) = (1-p)p^{k}, k=0,1,2,\dots
\]
Then they defined a growing process on $G$ according to the following rules: (1) at time zero, only the root (a fixed vertex of $G$)is declared infected, (2) at time $n \geq1$, a new vertex is infected if it is at graph distance at most $R_v$ of some vertex $v$ previously infected, and (3) infected vertices remain infected forever. They investigated the critical value $p_c(G)$ above which this process spreads indefinitely through the graph with positive probability.

They worked in a few settings including locally finite graphs in the sense that
\[ \Delta = \sup_{v \in {\mathcal G}} \{d(v)\}   < \infty \]
where $d(v)$ is the number of neighbors (or {\textit{degree}) of a vertex $v$.

An interesting question is whether such a model presents
\textit{phase transition} in the sense that for $p_c({\mathcal G}) := \inf\{p: \bbP(|I|=\infty\})$ we have that
\[0 < p_c({\mathcal G}) < 1.\]
They provided an answer which relies on a
comparison between the \dpm\ and the independent site 
percolation model. To understand this, consider $p_c^{site}({\mathcal G})$ 
the critical probability for the independent site percolation 
model on ${\mathcal G}$. 

\begin{teo}[Lebenstayn and Rodriguez~\cite{LR}]
\label{T: EP4}
Let ${\mathcal G}$ be of bounded degree ($\Delta < \infty$) and be such that $p_c^{site}({\mathcal G}) < 1$. Then
\begin{displaymath} 
0 < p_c({\mathcal G}) < 1.
\end{displaymath}
\end{teo}

The proof they presented relies on the following two propositions, the first
one is a comparison which gives an upper bound to $p_c({\mathcal G})$.

\begin{prop}[Lebenstayn and Rodriguez~\cite{LR}]
\label{P: EP2}
\begin{displaymath}
p_c({\mathcal G})  \leq p_c^{site}({\mathcal G})
\end{displaymath}
\end{prop}

\noindent
while the second one gives a lower bound for the case that 
${\mathcal G}$ is of bounded degree.

\begin{pro}[Lebenstayn and Rodriguez~\cite{LR}]
\label{P: EP3}
Suppose that ${\mathcal G}$ is a graph of bounded degree. Then
\begin{displaymath}
p_c({\mathcal G}) \ge -1 + {\left(1+ \frac{1}{\Delta -1}\right)}^{1/2}.
\end{displaymath}
\end{pro}

\subsection{Disk percolation on trees}
\label{S: Trees}\hfill

Consider a tree $\bbT$ (a connected graph with no cycles) and its set of vertices ${\mathcal V}(\bbT)$.
We say that a tree, $\bbT_d$, is \textit{homogeneous}, if each one of its vertices has degree (number of neighbours) $d+1$.
 
\begin{teo}[Lebenstayn and Rodriguez~\cite{LR}]
\label{T: EP10}
For any $d \ge 2$
\begin{displaymath}
-1 + \left(1-\frac{1}{d}\right)^{1/2} \leq p_c(\bbT_d) \leq 1 - \left(1-\frac{1}{d}\right)^{1/2}.
\end{displaymath}
\end{teo}

\begin{cor}[Lebenstayn and Rodriguez~\cite{LR}]
\label{C: EP11}
For any $d \ge 2$
\begin{displaymath}
p_c(\bbT_d) = {1/(2d)} + {\mathcal O}({1/{d^2}}) \text{ as } d \to \infty.
\end{displaymath}
\end{cor}

Single out one vertex from ${\mathcal V}(\bbT)$ and call this ${\mathcal O}$, the origin of ${\mathcal V}(\bbT)$. For each two vertices $u,v \in {\mathcal V}(\bbT)$, 
consider that $ u \leq v$ if $u$ belongs to the path connecting ${\mathcal O}$ to $v$. 

For a tree $\bbT$ and $n \geq 1$ we define

\begin{displaymath}
\bbT^u := \{ v \in \mathcal{V}: u \leq v \},
\end{displaymath}
\begin{displaymath}
\bbT_n^u := \{ v \in \bbT^u: d(v,{\mathcal O}) \leq d(u,{\mathcal O}) + n \}
\end{displaymath}
and
\begin{equation*}
\label{E: defPartialT}
M_n(u) := | \partial \bbT_n^u | := |\{ v \in \bbT^u: d(v,{\mathcal O}) = d(u,{\mathcal O}) + n\}|.
\end{equation*}

\begin{defn}
Let us define for a tree ${\mathbb T}$
\begin{displaymath}
{\textrm{dim\ inf\ }  \partial {\mathbb T}} : = \lim_{n \rightarrow \infty}
\min_{v \in \mathcal{V}} \frac{1}{n} \ln M_n(v).
\end{displaymath}
\end{defn}

Observe that 
\[ \textrm{dim\ inf\ } \partial {\mathbb T}_d = \ln d.\]

\begin{defn}
\label{D: SST}
We say that a tree, $\bbT_S$, is \textit{spherically symmetric}, if
any pair of vertices at the same distance from the origin, have the same
degree.
\end{defn}

\begin{teo}[Lebenstayn and Rodriguez~\cite{LR}]
\label{T: EP9}
For any spherically symmetric tree $\bbT_S$
\[ p_c (\bbT_S) \le 1-{\left(1- e^{-dim\ inf \partial \bbT_S}\right)}^{1/2} \]
\end{teo}

\section{Fireworks on $\bbN$}
\label{S: N}

The Fireworks processes are another interesting version of the 
\rpm. Junior \textit{et al}~\cite{Junior}  and 
Gallo \textit{et al}~\cite{GGJR} recently studied discrete time stochastic systems on $\bbN$ modeling 
processes of rumor spreading. In their models the involved individuals can either have an active 
role, working as spreaders and transmiting the information within a random distance to their right, 
or a passive role, hearing the information from spreaders within a random distance to their left. 
The appetite in spreading or hearing the rumor is represented by a set of random variables whose 
distributions may depend on the individuals positions on $\bbN$. Their main goal is to understand - 
based on the distribution of those random variables - whether the probability of having an infinite
set of individuals knowing the rumor is positive or not.

Junior \textit{et al}~\cite{Junior} manage to write the survival event as a limit of an increasing
sequence of events whose probability can be bounded by a nice use of FKG inequality. The
use of a non-standard version of Borel-Cantelli lemma helped in the task of finding conditions for the processes to die out. Gallo \textit{et al}~\cite{GGJR} based the proofs of their results on a clever relationship between the rumor
processes and a specific discrete time renewal process. With this technique they were able to obtain more precise results for homogeneous versions of the processes.

Consider $\{u_i\}_{i \in \bbN}$ a set of vertices of $\bbN$ such that
$ 0 < u_1 < u_2 <  \cdots $ and a set of independent
random variables $\{R_i\}_{i \in \bbN}$ assuming values in $\bbZ_+$.

\subsection{Fireworks}
\label{S: FP}

At time 0, information travels a distance $R_0$ towards the right side of the 
origin, in such
a way that all vertices $u_i \le R_0$ get informed. In general, at every
discrete time $t$ a vertex $u_j$ informed at time $t-1$ passes the information on (whithin $R_j$, its \textit{radius of influence}) and they do this just once, informing
the vertices $u_i$ (only those vertices which have not been informed
before) $u_j < u_i \le u_j + R_j.$ Observe that, except for the
set of vertices $\{u_i\},$ all other vertices are nonactionable, meaning that
their \textit{radius of influence} equals 0 almost surely.

\subsubsection{Homogeneous Fireworks}
\label{S: hoc}

Consider all the $R_i \sim R$ (having the same distribution) and $u_i=i$
for all $i$.

\begin{teo}[Junior \textit{et al}~\cite{Junior}]
\label{T: CriterioJunior}
Consider in the \hofp
\[a_n=\prod_{i=0}^{n}\bbP(R \leq i).\]

Then
\[ \sum_{n=1}^{\infty}a_n = \infty \hbox{ if and only if } \bbP[V]=0. \]
\end{teo}

\begin{teo}[Gallo \textit{et al}~\cite{GGJR}]
\label{T: CriterioGallo}
For the \hofp, 
\[ \bbP (V) = \left [ 1 + \sum_{j=1}^{\infty} \prod_{i=0}^{j-1}\bbP(R \leq i)  \right ]^{-1}.
\]
\end{teo}

Observe that the result presented in Theorem~\ref{T: CriterioJunior} is nicely generalized in Theorem~\ref{T: CriterioGallo}.

\begin{exa}
Consider the \hofp\ such that
\begin{equation*}
\label{E: Power Law}
\bbP(R = k) = \frac{2}{(k+2)(k+3)} \textrm {
for } k \in \bbN^*.
\end{equation*}
Then $\bbP[V]=\frac{1}{2}.$
\end{exa}

\begin{cor}[Junior \textit{et al}~\cite{Junior}]
\label{C: FPHo}
For the \hofp, consider
\[ L = \lim_{n \rightarrow \infty}n\bbP(R \geq n). \]

We have that
\begin{enumerate}
\item[\textit{(i)}] If $L > 1$ then $\bbP[V]>0.$
\item[\textit{(ii)}] If $L < 1$ then $\bbP[V]=0.$
\item[\textit{(iii)}] If $L = 1$ and there exists $N$ such that for all $n \geq N$
\begin{displaymath}
\bbP(R \geq n) \leq \frac{1}{n-1}, \textit{ then } \bbP[V]=0.
\end{displaymath}
\end{enumerate}
\end{cor}
Let $M$ be the final number of spreaders.
\begin{teo}[Gallo \textit{et al}~\cite{GGJR}]
\label{T: numinf}
If $\bbE(R) < \infty$ then the random variable $M$ has finite expectation. Besides, $M$ has exponential tail distribution when $\bbP(R \leq n)$ increases exponentially fast to $1$.
\end{teo}

Under more specific assumptions, it is possible to obtain more precise information on the tail distribution.
Items \textit{(i)} and \textit{(iii)} of next proposition follows from Proposition B.2 of Gallo \textit{et al.}~\cite{Gallo}, item \textit{(ii)} is due to Remark 5 from Bressaud \textit{et al.}~\cite{Bressaud} and item \textit{(iv)} follows from Theorem 1.1 of Garsia and Lamperti~\cite{Garsia}.

\begin{prop}[Gallo \textit{et al}~\cite{GGJR}]
\label{coro1}
 We have the following explicit bounds for the tail distributions.
\begin{enumerate}
\item[\textit{(i)}]  If $\bbP(R > k) \leq C_{r}r^{k},\,k\ge1$, for some $r\in(0,1)$ and a constant $C_{r}\in(0,\log \frac{1}{r})$ then
\[
\bbP(M \ge k)\leq \frac1{C_{r}}(e^{C_{r}}r)^{k}.
\]
\item[\textit{(ii)}] If $\bbP(R > k) \sim(\log k)^{\beta}k^{-\alpha}$, $\beta\in \bbR$, $\alpha>1$, then there exists  $C>0$ such that, for large $k$'s, we have
$\bbP(M\ge k)\leq C(\log k)^{\beta}k^{-\alpha}.$
\item[\textit{(iii)}] If $\bbP(R > k)= \frac{r}{k},\,k\ge1$ where $r\in(0,1)$, there exists  $C>0$ such that, for large $k$, we have
\[
\bbP(M\ge k)\leq C\frac{(\ln k)^{3+r}}{(k)^{2-(1+r)^{2}}}.
\]
\item[\textit{(iv)}] If $\bbP(R > k) \sim ((k+1)/(k+2))^{\alpha}$,  $\alpha\in(1/2,1)$, then there exists $C=C(\alpha)>0$ such that, for large $k$, we have
\[
\bbP(M\ge k)\leq \frac{C}{k^{1-\alpha}}.
\]
\end{enumerate}
\end{prop}

\subsubsection{Heterogeneous Fireworks}

\begin{rem}
\label{R: PFHeMorte} Consider the \hefp. One can get a
sufficient condition for $\bbP[V]=0$ ($\bbP[V]>0$) by a coupling argument.
Consider $\bbP(R_i \geq k ) \leq \mathbf{P}(R \geq k)$
($\bbP(R_i \geq k ) \geq \mathbf{P}(R \geq k)$)
for some random variable $R$ whose distribution $\mathbf{P}$ satisfies
$\lim_{n \rightarrow \infty}n \mathbf{P}(R \geq n) < 1$
($\lim_{n \rightarrow \infty}n \mathbf{P}(R \geq n) > 1$).
Finally use item ($ii$) (item ($i$)) of Corollary~\ref{C: FPHo}.
\end{rem}

\begin{teo}[Junior \textit{et al}~\cite{Junior}]
\label{T: PFHeVive}
Consider a \hefp\ for which actionable vertices are at integer positions
$u_0 = 0 < u_1 < u_2 < \dots $ such that $u_{n+1} - u_n \leq m$, for $m \geq 1.$
Besides, let us assume $\bbP (R_{n} < m) \in (0,1)$ for all $n.$
Then
\begin{enumerate}
\item[\textit{(i)}] If $\sum_{n=0}^{\infty}[\bbP(R_{n} < tm)]^t < \infty$ for some $t \geq 1$
then $\bbP[V]>0.$
\item[\textit{(ii)}] If for some random variable $R,$ with distribution $\mathbf{P},$
the following conditions hold
\begin{itemize}
\item $\mathbf{P}(R \geq k) - \bbP(R_{n} \geq k) \leq b_k$
for all $ k \geq 0$ and all $n \geq 0,$
\item $\lim_{n \rightarrow \infty}n[\mathbf{P}(R \geq n) - b_n] > m,$
\item $\lim_{n \rightarrow \infty}b_n = 0,$
\end{itemize}
then $\bbP[V]>0.$
\item[\textit{(iii)}] $\bbP(V) \geq
\prod_{j=0}^{\infty}\Big[1-\prod_{i=0}^{j}\bbP(R_{j-i} < (i+1)m)\Big].$
\end{enumerate}
\end{teo}

\subsection{Reverse Fireworks}
\label{S: RFP}

At time 0, only the origin has the information. At time 1, individuals placed at
vertices $u_i$ such that $u_i \le R_i $ get the information from the
origin. At time $t \in \bbN$
the set of vertices $ u_j $ which can find an informed individual at time $t-1$ within
a distance $R_j$ to its left, get the information. Let us call this set $A_t$.
If for some $t$, $A_t$ is empty the process stops. If the process never stops we say it survives.

Let $S$ be the event ``the reverse process survives". Besides, we denote by $Z$ the final number of spreaders. 

\subsubsection{Homogeneous Reverse Fireworks}
\label{S: horc}

Consider all the $R_i$ having the same distribution and $u_i=i$
for all $i$.

\begin{teo}[Junior {\it et al}~\cite{Junior}]
\label{T: horfpJunior}
Consider the \horfp. We have that
\begin{enumerate}
\item[\textit{(i)}] If $\mathbb{E}(R) = \infty$ then $\bbP(S) = 1.$
\item[\textit{(ii)}] If $\mathbb{E}(R) < \infty$ then $\bbP(S) = 0.$ 
\end{enumerate}
\end{teo}

\begin{teo}[Gallo {\it et al}~\cite{GGJR}]
\label{T: horfpGallo}
Consider the \horfp. If $\mathbb{E}(R) < \infty$ then
$Z \sim {\mathcal G}\displaystyle \left (\prod_{k=0}^{\infty}\bbP(R \leq k) \right )$ in the sense that for $p=\prod_{k=0}^{\infty}\bbP(R \leq k)$ we have 
\[\bbP(Z=k) = p (1-p)^k  \hbox{ for all } k.\]
\end{teo}

For any $n\ge1$, let $Z(n)$ be the number of spreaders in $\{1,\ldots,n\}$.
We will now state limit theorems for the proportion of spreaders within $\{1,\ldots,n\}$, $Z(n)/n$, when $n$ tends to $\infty.$

Let
\begin{align*}
\label{eq: mu1}
\mu&:=1+\sum_{j=1}^{\infty}\prod_{i=0}^{j-1}\bbP(R \leq i)\  \textrm{ and }\\
\sigma^2&:= \sum_{k=1}^{\infty}k^2\bbP(R > k-1)\prod_{i=0}^{k-2}\bbP(R \leq i)-\mu^2.
\end{align*}

Notice that $\mu < \infty$ implies that $ \displaystyle \prod_{k=0}^{\infty}\bbP(R \leq k)= 0$ (this implies $\bbE{(R)} = \infty$) .

\begin{teo}[Gallo {\it et al}~\cite{GGJR}]
\label{theo:Nn}If $\mu<\infty$ then
\[
\frac{Z(n)}{n}\stackrel{a.s.}{\longrightarrow}\mu^{-1},
\]
and thus, with probability one, $\mu^{-1}$ is the final proportion of spreaders.
Moreover, if $\sigma^2 \in(0,\infty)$,
then
\[
\sqrt{n}\left(\frac{Z(n)}{n}-\mu^{-1}\right)\stackrel{\mathcal{D}}{\rightarrow}\mathcal{N}\left(0,\frac{\sigma^2}{\mu^3}\right).
\]

Otherwise, $Z(n)/n \rightarrow 0$.
\end{teo}

In particular, observe that if the $\bbP(R \leq k)$'s satisfy at the same time $\prod_{k=0}^{\infty} \bbP(R \leq k)=0$ and $\mu=\infty$ (for instance, if they are as in items (iii) and (iv) of Proposition~\ref{coro1}), then the information  reaches infinitely many individuals, but the final proportion of informed individuals is zero.

\subsubsection{Heterogeneous Reverse Fireworks}
\begin{teo}[Junior {\it et al}~\cite{Junior}]
\label{T: PFRHe}
Consider the \herfp. It holds that
\begin{enumerate}
\item[\textit{(i)}] $\sum_{k=1}^{\infty}\bbP(R_{n+k} \geq k) = \infty$
for all $n$ if and only if $\bbP(S) = 1.$
\item[\textit{(ii)}] If $\sum_{n=1}^{\infty}\prod_{k=1}^{\infty}\bbP(R_{n+k} < k) < \infty$
then $\bbP(S)> 0.$
\end{enumerate}
\end{teo}

\begin{rem}
\label{R: PFRHe}
By a coupling argument and Theorem~\ref{T: horfpJunior} one can see that
if there is a random variable $R$, whose distribution is $\mathbf{P}$, with
$\bbE[R] < \infty$ ($\bbE[R] = \infty$), such that
$\bbP(R_{n} \geq k) \leq \mathbf{P}(R \geq k)$
($\bbP(R_{n} \geq k) \geq \mathbf{P}(R \geq k)$) for all $k$
then $\bbP(S) = 0$ ($\bbP(S) = 1$).
\end{rem}

\begin{exa}
\label{E: um}
It is possible to have in the \hefp\  the expectation of the \textit{radius of
influence} infinite for all vertices toghether and the process dies out almost surely. 

Let $\{b_n\}_{n \in \bbN}$ be a non-increasing sequence convergent to 0 and such that $b_0<1.$
\begin{itemize}
\item[\textit{(i)}] $\bbP(R_n = 0) = 1 - b_n$ and $\bbP(R_n = k) = b_{n+k-1} - b_{n+k}$ for $k \geq 1.$
\item[\textit{(ii)}] $\sum_{n=0}^{\infty}b_n = \infty.$
\item[\textit{(iii)}] $\lim_{n \rightarrow \infty} n b_n = 0.$
\end{itemize}
Observe that $\mathbb{E}(R_n) = \infty$ for all $n$ from \textit{(ii)}. Besides
$\bbP[V]=0$ from \textit{(iii)}, because
For
\[V_n = \{ \hbox{The individual at vertex } u_n \hbox{ gets the information}\},\] 
\begin{equation}
\label{E: example}
\bbP(V_n) \leq \sum_{k=0}^{n-1}\bbP(R_k \geq n-k) =
\sum_{k=0}^{n-1}b_{n-1} = (n-1)b_n.
\end{equation}
and the fact that $V = \lim_{n \to \infty} V_n.$
\end{exa}
\begin{exa}
It is possible to have in the \hefp\ the expectation of the \textit{radius of
influence} finite for all vertices and the process survives with positive probability.
Assume that $\sum_{n=0}^{\infty}b_n < \infty,$ while
\begin{itemize}
\item[\textit{(i)}] $\bbP(R_n = 0) = b_n$
\item[\textit{(ii)}] $\bbP(R_n = 1) =1 - b_n$
\end{itemize}
Then $\mathbb{E}(R_n) < 1$ for all $n$ and
$\mathbb{P}(V) > 0$ by item ($i$) of Theorem~\ref{T: PFHeVive}
with $m=t=1$.
\end{exa}
\begin{exa} Next we present an example where $\bbP[S]=1$ for
a \herfp\ while $\bbP[V]=0$ for a \hefp. For this aim consider
\begin{enumerate}
\item[\textit{(i)}] $\bbP(R_n = 0) = 1 - b_n$ and $\bbP(R_n = n) = b_n.$
\item[\textit{(ii)}] $\sum_{n=0}^{\infty}b_n = \infty.$
\item[\textit{(iii)}] $\lim_{n \rightarrow \infty} n b_n = 0.$
\end{enumerate}
Observe that even though $\lim_{n \to \infty} \bbE[R_n] = 0$ and
$\lim_{n \to \infty} \bbP(R_n = 0) = 1,$
from Theorem~\ref{T: PFRHe} and \textit{(ii)} it is true for the \herfp\
that $\mathbb{P}(S) = 1.$ In the opposite direction, by~(\ref{E: example})
and \textit{(iii)} one have that $\bbP[V]=0$ for the \hefp.
\end{exa}

%%%%%%%%%%%
\section{Cone percolation on $\bbT_d$}
\label{S: MR}

Junior \textit{et al}~\cite{Junior2} consider a process which allows us to associate the dynamic activation on the set of vertices to a discrete rumor process. Individuals become spreaders as soon as they hear the rumor. Next time, they propagate the rumor within their \textit{radius of influence} and immediately become stiflers. Junior \textit{et al}~\cite{Junior2} establish whether the process has positive probability of involving an infinite set of individuals. Besides, they present sharp lower and upper bounds for the probability of that event, depending on the general distribution of the random variables that define the \textit{radius of influence} of each individual. Their
proofs are based on comparisons with branching processes.

Pick a $v \in {\mathcal V}(\bbT_d)$ such that $d({\mathcal O},v)=1$ and consider 
$\bbT^+_d = \bbT_d \backslash \bbT_d^+(v).$ 
Consider $\bbP_+$ and $\bbP$ the probability measures associated to the processes on $\bbT_d^+$ and $\bbT_d$
(we do not mention the random variable $R$ unless absolutely necessary). By
a coupling argument one can see that for a fixed distribution of $R$
\begin{equation*}
\bbP_+[V] \leq \bbP[V].
\end{equation*}

Furthermore, by the definition of $\bbT_d^+$ and its relation with $\bbT_d$ we have that
for a fixed distribution of $R$
\begin{equation*}
\label{E: Equal}
\bbP_+[V]=0 \hbox{ if and only if } \bbP[V]=0.
\end{equation*}

Let $p_0 = \bbP(R=0).$

\begin{teo}[Junior {\it et al}~\cite{Junior2}]
\label{T:CSPAH}
Consider the \cpp\ on $\bbT_d^+$ with \textit{radius of influence} $R.$ 
\begin{enumerate}
\item[\textit{(i)}] If $(1-p_0) d > 1$, then $\bbP_+[V] > 0,$
\item[\textit{(ii)}] If $(1-p_0) d \le 1$ and $\bbE(d^R) > 1 + p_0$, then $\bbP_+[V] > 0,$
\item[\textit{(iii)}] If $\bbE(d^R) \leq 2 - \frac{1}{d}$, then $\bbP_+[V] = 0.$
\end{enumerate}
\end{teo}

\begin{teo}[Junior {\it et al}~\cite{Junior3}]
\label{T: THT1}
Consider a \cpp\  on $\bbT_d$.
Then for $\mathbb{E}(d^R) < 2 - \frac{1}{d}$, we have
\begin{displaymath}
 \frac{d + \mathbb{E}\displaystyle \left (d^R \right) - p_0 }{d[1-\mathbb{E}\displaystyle \left (d^R \right) + p_0]} \leq \mathbb{E}(|I|)  \leq \frac{\mathbb{E}\displaystyle \left (d^R \right)+d-2}{2d - 1 - d\mathbb{E}\displaystyle \left (d^R \right)}.
\end{displaymath}
\end{teo}

\begin{exa}[Junior {\it et al}~\cite{Junior3}] Consider $R \sim {\mathcal G}(1-p)$, a \textit{radius of
influence} satisfying
\[
\bbP(R = k) = (1-p)p^{k}, k=0,1,2,\dots
\]
and assume also $pd < \frac{1}{2}$.
So we have
\begin{displaymath}
\frac{1-dp +p -p^2}{1-2dp+dp^2} \leq \mathbb{E}(|I|) \leq \frac{1-dp-p}{1-2dp}.
\end{displaymath}
That gives us a fairly sharp bound even when we pick $p$ and $d$ such that $pd$ is very close to $\frac{1}{2}$ as, for example, $p = 10^{-6}$ and $d = 499,000$. 
For these parameters we get $250.438 \leq \mathbb{E}(|I|) \leq 250.501$.
\end{exa}

Let $\rho$ and $\psi$ be, respectively, the smallest non-negative roots of the equations
\begin{align*}
& \bbE(\rho^{d^R}) + (1 - \rho)p_0 = \rho, \\
& \bbE(\psi^{\frac{d}{d-1}(d^{R}-1)}) = \psi.
\end{align*}

\begin{teo}[Junior {\it et al}~\cite{Junior2}]
\label{T: SobrevivenciaTd+}
Consider the \cpp\ on $\bbT_d^+.$
Then
\begin{displaymath}
1 - \rho \leq \bbP_+(V) \leq 1 - \psi.
\end{displaymath}
\end{teo}

\begin{teo}[Junior {\it et al}~\cite{Junior2}] 
\label{T: ViveTd}
For the \cpp\ on $\bbT_d$ with \textit{radius of influence} $R$, it holds that
\begin{equation*}
\label{E: ViveTd}
1 - \displaystyle \left(1 -
\rho^{\frac{d+1}{d}}\right)p_0 -
\bbE\displaystyle \left(\rho^{\frac{(d+1)}{d}d^{R}}\right)
\leq \bbP[V] \leq 1 - \bbE\displaystyle
\left(\psi^{\frac{(d+1)}{d-1}(d^{R}-1)}\right).
\end{equation*}
\end{teo}

Consider $d=2$ and $R$ following a Binomial distribution with parameters
$4$ and $\frac{1}{2}$
($R \sim \mathcal{B}(4,\frac{1}{2})$).
Therefore $\rho$ and $\psi$ are, respectively,
solutions of
\begin{align*}
& x^{16} + 4x^8 + 6x^4 + 4x^2 - 16x + 1 = 0, \\& x^{30} + 4x^{14} +
6x^6 + 4x^2 - 16x + 1 = 0.
\end{align*}
So $\rho$ = 0.0635146 and $\psi$ = 0.06350850, which implies that
\[
0.937435919 \leq \bbP[V] \leq 0.937435962.
\]

\begin{teo}[Junior {\it et al}~\cite{Junior2}] 
\label{aest}
The \hcpp\ in $\bbT^+_d$ has a giant component with positive pro\-ba\-bi\-li\-ty if
for some fixed $n$,
\begin{equation}
\label{eq1}
\liminf_{j \rightarrow \infty} d^n
\prod_{k=0}^{n-1}[1-\prod_{i=0}^{k}\bbP_+[R_{jn+i} < k+1-i]] >
1.
\end{equation}
\end{teo}

A consequence of Theorem 5.4 from Bertacchi and Zucca~\cite{BZ} is the following result

\begin{cor} 
\label{C: novo}
Consider a \horfp\ on $\bbT_d$. Then
\[ \bbP(S) = 1 \textrm { if and only if } \sum_{n=1}^{\infty} d^n \bbP(R \geq n) = \infty.
\]
\[ \bbP(S) = 0 \textrm { if and only if } \sum_{n=1}^{\infty} d^n \bbP(R \geq n)\prod_{j=1}^{n-1}[1 - \bbP(R \geq j)] \leq 1
\]
\end{cor}
%%%%%%%%%%%

\begin{teo}[Junior {\it et al}~\cite{Junior3}]
\label{T: esferic}
For a \cpp\ in $\bbT_S$ and $R$, the \textit{radius of influence}, $\bbP(V) > 0$ if
\begin{equation*}
\label{E: esferic}
\lim_{n \rightarrow \infty}\sqrt[n]{\rho_n} > 
e^{-\textrm{dim\ inf\ } \partial {\mathbb T}_S}
\end{equation*}
where
\begin{equation*}
\rho_n := \prod_{k=0}^{n-1} [1-\prod_{i=0}^{k}\bbP(R < i+1)].
\end{equation*}
\end{teo}

\begin{cor}[Junior {\it et al}~\cite{Junior3}]
\label{C: esferic1}
For a \cpp\ in $\bbT_S$ and $R$, a \textit{radius of influence} satisfying $\bbP(R \leq k) = 1$ for some $k \in \bbN$, $\bbP(V) > 0$ if

\begin{equation*}
\label{E: esferic1}
\textrm{dim\ inf\ } \partial {\mathbb T}_S > \ln \displaystyle \left[\frac{1}{1 - \prod_{j=1}^{k}\bbP(R < j)} \right].
\end{equation*}
\end{cor}

\begin{defn}
A $k$\textit{-periodic tree} with degree $\tilde{d} = (d_1, \cdots, d_k)$, $d_i \geq 2$ for all 
$i=1,2,\cdots, k$, is as tree such that for any vertex whose 
distance to the origin is $nk+i-1$ for some $ n \in\bbN $ has degree $d_i + 1$. 
We refer to this tree as $\bbT_{\tilde{d}}$.
\end{defn}

\begin{exa}[Junior {\it et al}~\cite{Junior3}]
\label{E: B(p)}
Consider  a \cpp\ in $\bbT_S$ with $R \sim {\mathcal B}(p),$ a \textit{radius of influence} sa\-tis\-fying
\[
\bbP(R = 1) = p = 1 - \bbP(R = 0).
\]
\begin{itemize}
\item[\textit{(i)}] If dim\ inf\ $\partial {\mathbb T}_S > -\ln p$ then $\bbP(v)>0,$
%\item If $\bbT_S = \bbT_{\tilde{d}}$ and $ G(\tilde{d}) > \frac{1}{p}$ then $\bbP(v)>0.$
\item[\textit{(ii)}] If $\bbT_S = \bbT_{\tilde{d}}$ and $ \sqrt[k]{\prod_{j=1}^{k}d_{j}} > \frac{1}{p}$ then $\bbP(V)>0.$
\end{itemize}
\end{exa}

\section{Random environments}
\label{S: RPRE}

In this section we review the Fireworks and the Reverse Fireworks processes, with a random number of stations at each vertex. Bertacchi and Zucca~\cite{BZ} consider an extra source
of randomness: the number of individuals sitting on each vertex. They consider
two families of random variables $\{N_x\}_{x \in {\mathcal G}}$ and 
$\{R_{x,i}\}^{i \in \bbN}_{x \in {\mathcal G}}$ such that $\{N_x, R_{x,i}\}$
are independent and $\{R_{x,i}\}_{i \in {\bbN}}$ are identically distributed for all
$x \in {\mathcal G}$ that is $R_{x,i} \sim R_x.$ In their paper $N_x$ represents
the random number of individuals at vertex $x$ (in particular
$N_{\mathcal O}$ is the number of individuals at the origin) while $\{R_{x,i}\}_{i=1}^{N_x}$
are their \textit{radius of influence}. The main question about this model is to understand under which conditions, the signal, starting from one vertex of a graph ($\bbN$ or a Galton-Watson tree), will spread indefinitely with positive probability or die out almost surely in a finite number of steps.

Bertacchi and Zucca~\cite{BZ} rely in their analisys on associating the processes with random numbers of stations (fireworks or 
reverse fireworks), with processes with one station per vertex as in Junior \textit{et al}~\cite{Junior} . Indeed, they consider processes with one station on each vertex $x$ and \textit{radius of influence}
${\tilde{R}}_x = {\mathbf 1}_{\{N_x \ge 1\}} \max\{R_{x,j}: j = 1, \dots, N_x\}.$ 
They call this process, the \textit{deterministic counterpart} or \textit{annealed counterpart} of the
original process. They observe that the \textit{annealed counterpart} does not retain any information about the environment,
nevertheless the probability of survival for the original process and for its \textit{annealed counterpart} are the same.

\subsection{Fireworks}

For $x \in {\mathcal G},$ let us define
\begin{displaymath}
\varphi_{N_x}(t) := \bbE(t^{N_x}) = \sum_{j=0}^{\infty} \bbP(N_x=j) t^j
\end{displaymath}

\subsubsection{Homogeneus Fireworks}
Consider $R_i \sim R$ and $N_x \sim N$ for all $x \in {\mathcal G}.$
Let us define 
\begin{displaymath}
f_{R,N}(n) := n \{1 - \varphi_{N}(\mathbb{P}(R < n))\}.
\end{displaymath}

\begin{teo}[Bertacchi and Zucca~\cite{BZ}]
\label{T: NAT} \hfill
\begin{itemize}
\item[\textit{(i)}]
If $ \displaystyle \limsup_{n \rightarrow \infty}f_{R, N}(n) < 1$
 then $\mathbb{P}(V) = 0$.
\item[\textit{(ii)}]
If $ \displaystyle \liminf_{n \rightarrow \infty}f_{R, N}(n) > 1$
 then $\mathbb{P}(V) > 0$.
\item[\textit{(iii)}] If $\mathbb{E}(N) < \infty$ and $\displaystyle \limsup_{n \to \infty} n \mathbb{P}( R \geq n) < \frac{1}{\mathbb{E}(N)}$ then $\mathbb{P}(V) = 0$.
\item[\textit{(iv)}] If $\mathbb{E} (N) < \infty$  and $\mathbb{E} (R) < \infty$ then $\mathbb{P}(V)  = 0$.
\item[\textit{(v)}] If $ \displaystyle \liminf_{n \rightarrow \infty}n \mathbb{P}( R \geq n)\varphi^{\prime}_{N}(\mathbb{P}(R < n)) > 1$
 then $\mathbb{P}(V) > 0$.
\end{itemize}
\end{teo}

A consequence of Theorem 1 from Gallo \textit{et al}~\cite{GGJR} is the following result

\begin{cor}
\begin{displaymath}
\mathbb{P}(V) = \left [ 1 + \sum_{j=1}^{\infty} \prod_{i=0}^{j-1}\displaystyle \varphi_{N} \left ( \displaystyle \mathbb{P}(R \leq
i) \right )\right ]^{-1}
\end{displaymath}
\end{cor}

%\begin{rem}
%\label{R: BZ-RandN} \hfill
%%\begin{cor} [Bertacchi and Zucca~\cite{BZ}] \hfill
%Proposition~\cite{P: NAT} admitss a similar corolary
%\begin{enumerate}
%\item
%For every unbounded random variable $R$ there exists a
%random variable $N$ such that $\mathbb{P}(v) > 0$.
%For $\epsilon > 0$ and $\delta \in (0,1)$ consider $N$ satisfying
%\[ \mathbb{P} \left (N \geq \frac{\ln(1-\delta)}{\ln (\mathbb{P} (R < n))} \right) \geq \frac{1 + \epsilon}{n \delta.}
%\]
%\item For every random variable $N$ such that $\mathbb{P}(N = 0) < 1$ there exists a random variable $R$ such that $\mathbb{P}(v) > 0$. Take $R$ satisfying $\mathbb{P} (R \geq n) = p_n$, where $p_n = \inf \{ t \geq 0; \varphi_{N} \left ( 1-t\right ) \leq 1 - \frac{2}{n} $\}
%\end{enumerate}
%%\end{cor}
%\end{rem}

\begin{rem}
It is possible to have $\mathbb{E}(N) = \infty$, $\mathbb{E}(R) = \infty$ and $\mathbb{P}(V) = 0$.  Take $\mathbb{P}( N \geq n) \sim \frac{1}{n}$ when $n \to \infty$ and $\mathbb{P}( R \geq n)  = \frac{1}{n \ln n \ln (ln n)}$
\end{rem}

\subsubsection{Heterogeneous Fireworks}

\begin{teo}[Bertacchi and Zucca~\cite{BZ}]
In the heterogeneous case, if
\[ \sum_{n=0}^{\infty} \prod_{i=0}^{n} \varphi_{N_i} \left ( \displaystyle \mathbb{P}(R_i <
n-i+1) \right )
\]
then $\mathbb{P}(V) > 0$.
\end{teo}

Adapting the proof of Theorem 2.3 from Junior {\it et al}~\cite{Junior} we have

\begin{teo}
In the heterogeneous case, if
\begin{itemize}
\item[\textit{(i)}] $ \varphi_{N_i}(\mathbb{P}(R_i < 1)) \in (0,1).$
\item[\textit{(ii)}] $ \lim_{n \to \infty} \prod_{i=0}^{n-1} \varphi_{N_i} \left ( \displaystyle \mathbb{P}(R_i <
2n-1) \right ) = 1.$
\item[\textit{(iii)}] $ \lim_{n \to \infty} \prod_{i=n}^{2n-1} \varphi_{N_i} \left ( \displaystyle \mathbb{P}(R_i <
2n-1) \right ) > 0.$
\end{itemize}
then $\mathbb{P}(V) = 0$.
\end{teo}

\subsection{Reverse Fireworks}

\subsubsection{Homogeneous Reverse Fireworks}

Let us define 

\[ W = \sum_{n=0}^{\infty} \left [ 1 - \varphi_{N} \left ( \displaystyle \mathbb{P}(R <
n) \right )   \right ]
\]

\begin{teo}[Bertacchi and Zucca~\cite{BZ}]
\label{T: NATR}\hfill
\begin{itemize}
\item[\textit{(i)}]
If $W = \infty$ then $\mathbb{P}(S) = 1.$
\item[\textit{(ii)}]
If $W < \infty$ then $\mathbb{P}(S) = 0.$
\end{itemize}
\end{teo}

Theorem~\ref{T: NATR} can also be obtained as a consequence of Theorem 3.2 from Junior {\it et al}~\cite{Junior2} or as a consequence of Theorem 2 from Gallo {\it et al}~\cite{GGJR}.

\begin{rem}[Bertacchi and Zucca~\cite{BZ}]
\label{R: BZ-RandN} 
Theorems~\ref{T: NAT} and~\ref{T: NATR} admit a similar corolary
\begin{itemize}
\item[\textit{(i)}]
For every unbounded random variable $R$ there exists a
random variable $N$ such that $\mathbb{P}(V) > 0$  ($\mathbb{P}(S) = 1$).
For $\epsilon > 0$ and $\delta \in (0,1)$ consider $N$ satisfying
\[ \mathbb{P} \left (N \geq \frac{\ln(1-\delta)}{\ln (\mathbb{P} (R < n))} \right) \geq \frac{1 + \epsilon}{n \delta}.
\]
\item[\textit{(ii)}] 
For every random variable $N$ such that $\mathbb{P}(N = 0) < 1$ there exists a random variable $R$ such that $\mathbb{P}(V) > 0$ ($\mathbb{P}(S) = 1$). Take $R$ satisfying $\mathbb{P} (R \geq n) = p_n$, where $p_n = \inf \{ t \geq 0; \varphi_{N} \left ( 1-t\right ) \leq 1 - \frac{2}{n} $\}.
\end{itemize}
%\end{cor}
\end{rem}

\subsubsection{Heterogeneous Reverse Fireworks}

\begin{teo}[Bertacchi and Zucca~\cite{BZ}]
In the heterogeneous case,
\[ \sum_{k=0}^{\infty} \left [ 1 - \varphi_{N_{n+k}} \left ( \displaystyle \mathbb{P}(R_{n+k} <
k) \right ) \right ] = \infty, \textrm { if and only if }  \mathbb{P}(S) = 1.
\]
By other hand, if
\[ \sum_{n=0}^{\infty} \prod_{k=1}^{\infty} \varphi_{N_{n+k}} \left ( \displaystyle \mathbb{P}(R_{n+k} <
k) \right ) < \infty, \mathbb{P}(S) > 0.
\]

\end{teo}

\subsection{Galton Watson}
Let us define the space of unlabelled GW-trees (the usual GW-trees). Consider a GW-process,
with offspring distribution $ \bbP(D = d),\ 0 \le d < \infty $. We assume that $\bbP(D = 1) < 1$ (otherwise the
resulting random tree is $\bbN$) and we suppose that $ \mu_D := \displaystyle \sum_{d=0}^{\infty} d \bbP(D = d) > 1$ (the supercritical case).
The underlying random graph will be a GW-tree generated by this process. Let $ \displaystyle g(s): = \sum_{d=0}^{\infty} s^d \bbP(D = d) $ be the generating function of $D$ and let $\pi \in [0, 1]$ be the smallest nonnegative fixed point of $g$. If $\bbP(D > k) = 0$ for some $k$ we say that the GW-tree has maximum degree $k$ or that it is $k$-bounded.

\subsubsection{Homogeneous Fireworks}
In this case, the random number of stations are independent and identically distributed $\bbN$-valued random variables
with common law $N$. Analogously, The radii of the stations are independent and identically distributed with distribution $R$ (either discrete or continuous random variable).

\begin{defn}
We define \[ \Phi(t) : = \varphi_N(\bbP(R < 1)) + \sum_{n=1}^{\infty}[\varphi_N(\bbP(R < n+1))-\varphi_N(\bbP(R < n))]t^n.\]
\end{defn}

In particular observe that
\[ \Phi(0) = \varphi_N(\bbP(R < 1)) \textrm {  and the case } N =1 \textrm { a.s.},\]
\[ \Phi(t) = \sum_{n=0}^{\infty}[\bbP( n \leq R < n+1) ]t^n.   \]

\begin{teo}[Bertacchi and Zucca~\cite{BZ}]
Consider a \hofp. We have that
\begin{itemize}
\item[\textit{(i)}] If $\Phi(\mu_D)- 1 > \Phi(0) = \varphi_N(P(R < 1))$ and $\bbP(N = 0) = 0 $ then for the Fireworks process there
is survival with positive probability for almost every realization of the environment such that
the underlying tree is infinite and there is at least one station at the root.
\item[\textit{(ii)}] If $\Phi(\mu_D) - 1 > \Phi(0) = \varphi_N(P(R < 1))$ and $\bbP(N = 0) > 0$ then for the 
Fireworks process
$\bbP(V|\tau  = T,N_{\mathcal{O}} = n) > 0$ for almost every $(T, n)$ such that $T$ is an infinite
(unlabelled) tree and $n \geq 1$.
\item[\textit{(iii)}] If the GW-tree is $k$-bounded and $\Phi(k) \leq 2- \frac{1}{k}$ then the Fireworks process becomes
extinct a.s. for almost every realization of the environment.
\end{itemize}
\end{teo}

\subsubsection{Homogeneous Reverse Fireworks}
In this case, the random number of stations are independent and identically distributed $\bbN$-valued random variables
with common law $N$, except by numbers of station the root $\mathcal{O}$. For the root, we take $N_{\mathcal{O}} = \min\{ n > 0: \bbP(N = n) > 0 \} $ Besides, the radii of the stations are independent and identically distributed with distribution $R$ (either discrete or continuous random variable).

\begin{defn}
We define \[ \phi_1(t) : = \sum_{n=1}^{\infty}[1 - \varphi_N(\bbP(R < n)) ]{\mu_D}^n \]
\[ \phi_2(t) : = \sum_{n=1}^{\infty}[1 - \varphi_N(\bbP(R < n)) ]{\mu_D}^n \prod_{j=1}^{n-1}\varphi_N(\bbP(R < j))  \]
\end{defn}

\begin{teo}[Bertacchi and Zucca~\cite{BZ}]
Consider a \horfp. The following hold
\begin{itemize}
\item[\textit{(i)}] If $\phi_1(\mu_D) = \infty$ then there is survival with probability 1 for the Reverse Fireworks process for
almost all realizations of the environment such that the underlying tree is infinite.
\item[\textit{(ii)}] If $\bbP(N = 0) = 0$, $\phi_1(\mu_D) < \infty$ and $\phi_2(\mu_D) > 1$ then there is survival with positive probability
(strictly smaller than 1) for the Reverse Fireworks process for almost all realizations of the
environment such that the underlying tree is infinite.
\item[\textit{(iii)}] If $\bbP(N = 0) > 0$, $\phi_1(\mu_D) < \infty$ and $\phi_2(\mu_D) > 1$ then $\bbP(S|\tau = T) \in (0, 1)$ for almost
every infinite (unlabelled) tree $T$.
\item[\textit{(iv)}] If $\phi_1(\mu_D) < \infty$ and $\phi_2(\mu_D) \leq 1$ then there is a.s. extinction for the Reverse Fireworks process
for almost all realizations of the environment.
\end{itemize}
\end{teo}

\begin{defn}
We define \[ \mathcal{M}_c := \left [ \limsup_{n \to \infty} \sqrt[n]{1 - \varphi_N(\bbP(R < n))}   \right ]^{-1}.\]
\end{defn}

\begin{cor}[Bertacchi and Zucca~\cite{BZ}]
There exists a critical value $\mu_c \in [1,\infty), \mu_c \leq  \mathcal{M}_c$ such that
\begin{itemize}
\item[\textit{(i)}] $\mu_D < \mu_c$ implies a.s. extinction for almost all realizations of the environment.
\item[\textit{(ii)}] $\mu_c < \mu_D < \mathcal{M}_c$ and $\bbP(N = 0) = 0$ implies survival with positive probability for almost all
realizations of the environment such that the underlying tree is infinite.
\item[\textit{(iii)}] $\mu_c < \mu_D < \mathcal{M}_c$ and $\bbP(N = 0) > 0$  implies survival with positive probability for almost every
infinite (unlabelled) tree.
\item[\textit{(iv)}] $\mathcal{M}_c < \mu_D$ implies survival with probability 1 for almost all realizations of the environment
such that the underlying tree is infinite.
\item[\textit{(v)}] If $\mu_D = \mu_c < \mathcal{M}_c$ then there is a.s. extinction for almost all realizations of the environment.
\end{itemize}
\end{cor}

\section{Open problems}
\label{S: OQ}

Some natural extensions for these models are those considering

\begin{itemize}
\item[\textit{(i)}] Fireworks processes (direct and reverse) on ${\bbZ^d}^+$. An especially interesting case is when $d=2$ and the \textit{boxes of influence} are distributed as $[0, R_x) \times [0, R_y)$ with $R_x$ independent of $R_y$ and the rumor starting from $(0,0)$ or from
every $(x,y)$ such that $x=0$ or $y=0$;
\item[\textit{(ii)}] Fireworks processes on $\bbZ$. Heterogenous versions with \textit{radius of influence} non i.i.d. and with in\-di\-vi\-duals being initially placed following a renewal process or a Markovian process.
\item[\textit{(iii)}] Reverse fireworks processes on $\bbZ$. Individuals throw their \textit{radius of influence} to every direction as in~(\ref{E: defBu}) (See Gallo \textit{et all}~\cite{GGJR}). They believe that conditions for survival will be the same but the final proportion of informed
individual will be strictly larger.
\item[\textit{(iv)}] Cone Percolation on Spherically Symmetric and on Galton Watson trees. Lower
and upper bounds for the survival probability and for the extinction time.
\end{itemize}

\noindent
Acknowledgments: F.P.M. wishes to thank NYU-Shanghai China and V.V.J. and K.R. wish to 
thank Instituto de Matem\'atica e  Estat\'{\i}stica-USP Brazil for kind hospitality.

\end{document}